\documentclass[10pt,a4paper]{article}
\usepackage{amsmath, amssymb}
\usepackage{latexsym, color}
\usepackage{epsfig}

\numberwithin{equation}{section}

\begin{document}

\newtheorem{thm}{Theorem}[section]
\newtheorem{cor}[thm]{Corollary}
\newtheorem{lem}[thm]{Lemma}
\newtheorem{prop}[thm]{Proposition}
\newtheorem{definition}[thm]{Definition}
\newtheorem{rem}[thm]{Remark}
\newtheorem{Ex}[thm]{EXAMPLE}
\newtheorem{Nota}[thm]{Notation}
\def\nm{\noalign{\medskip}}

\bibliographystyle{plain}


\newcommand{\qed}{\hfill \ensuremath{\square}}
\newcommand{\ds}{\displaystyle}
\newcommand{\pf}{\noindent {\sl Proof}. \ }
\newcommand{\p}{\partial}
\newcommand{\pd}[2]{\frac {\p #1}{\p #2}}
\newcommand{\norm}[1]{\| #1 \|}
\newcommand{\dbar}{\overline \p}
\newcommand{\eqnref}[1]{(\ref {#1})}
\newcommand{\na}{\nabla}
\newcommand{\ep}{\epsilon}
\newcommand{\one}[1]{#1^{(1)}}
\newcommand{\two}[1]{#1^{(2)}}

\newcommand{\Abb}{\mathbb{A}}
\newcommand{\Cbb}{\mathbb{C}}
\newcommand{\Ibb}{\mathbb{I}}
\newcommand{\Nbb}{\mathbb{N}}
\newcommand{\Kbb}{\mathbb{K}}
\newcommand{\Rbb}{\mathbb{R}}
\newcommand{\Sbb}{\mathbb{S}}

\renewcommand{\div}{\mbox{div}~}

\newcommand{\la}{\langle}
\newcommand{\ra}{\rangle}

\newcommand{\Hcal}{\mathcal{H}}
\newcommand{\Lcal}{\mathcal{L}}
\newcommand{\Kcal}{\mathcal{K}}
\newcommand{\Dcal}{\mathcal{D}}
\newcommand{\Pcal}{\mathcal{P}}
\newcommand{\Qcal}{\mathcal{Q}}
\newcommand{\Scal}{\mathcal{S}}

\def\Ba{{\bf a}}
\def\Bb{{\bf b}}
\def\Bc{{\bf c}}
\def\Bd{{\bf d}}
\def\Be{{\bf e}}
\def\Bf{{\bf f}}
\def\Bg{{\bf g}}
\def\Bh{{\bf h}}
\def\Bi{{\bf i}}
\def\Bj{{\bf j}}
\def\Bk{{\bf k}}
\def\Bl{{\bf l}}
\def\Bm{{\bf m}}
\def\Bn{{\bf n}}
\def\Bo{{\bf o}}
\def\Bp{{\bf p}}
\def\Bq{{\bf q}}
\def\Br{{\bf r}}
\def\Bs{{\bf s}}
\def\Bt{{\bf t}}
\def\Bu{{\bf u}}
\def\Bv{{\bf v}}
\def\Bw{{\bf w}}
\def\Bx{{\bf x}}
\def\By{{\bf y}}
\def\Bz{{\bf z}}
\def\BA{{\bf A}}
\def\BB{{\bf B}}
\def\BC{{\bf C}}
\def\BD{{\bf D}}
\def\BE{{\bf E}}
\def\BF{{\bf F}}
\def\BG{{\bf G}}
\def\BH{{\bf H}}
\def\BI{{\bf I}}
\def\BJ{{\bf J}}
\def\BK{{\bf K}}
\def\BL{{\bf L}}
\def\BM{{\bf M}}
\def\BN{{\bf N}}
\def\BO{{\bf O}}
\def\BP{{\bf P}}
\def\BQ{{\bf Q}}
\def\BR{{\bf R}}
\def\BS{{\bf S}}
\def\BT{{\bf T}}
\def\BU{{\bf U}}
\def\BV{{\bf V}}
\def\BW{{\bf W}}
\def\BX{{\bf X}}
\def\BY{{\bf Y}}
\def\BZ{{\bf Z}}


\newcommand{\Ga}{\alpha}
\newcommand{\Gb}{\beta}
\newcommand{\Gd}{\delta}
\newcommand{\Ge}{\epsilon}
\newcommand{\Gve}{\varepsilon}
\newcommand{\Gf}{\phi}
\newcommand{\Gvf}{\varphi}
\newcommand{\Gg}{\gamma}
\newcommand{\Gc}{\chi}
\newcommand{\Gi}{\iota}
\newcommand{\Gk}{\kappa}
\newcommand{\Gvk}{\varkappa}
\newcommand{\Gl}{\lambda}
\newcommand{\Gn}{\eta}
\newcommand{\Gm}{\mu}
\newcommand{\Gv}{\nu}
\newcommand{\Gp}{\pi}
\newcommand{\Gt}{\theta}
\newcommand{\Gvt}{\vartheta}
\newcommand{\Gr}{\rho}
\newcommand{\Gvr}{\varrho}
\newcommand{\Gs}{\sigma}
\newcommand{\Gvs}{\varsigma}
\newcommand{\Gj}{\tau}
\newcommand{\Gu}{\upsilon}
\newcommand{\Go}{\omega}
\newcommand{\Gx}{\xi}
\newcommand{\Gy}{\psi}
\newcommand{\Gz}{\zeta}
\newcommand{\GD}{\Delta}
\newcommand{\GF}{\Phi}
\newcommand{\GG}{\Gamma}
\newcommand{\GL}{\Lambda}
\newcommand{\GP}{\Pi}
\newcommand{\GT}{\Theta}
\newcommand{\GS}{\Sigma}
\newcommand{\GU}{\Upsilon}
\newcommand{\GO}{\Omega}
\newcommand{\GX}{\Xi}
\newcommand{\GY}{\Psi}

\newcommand{\beq}{\begin{equation}}
\newcommand{\eeq}{\end{equation}}

\title{Characterization of the electric field concentration between two adjacent spherical perfect conductors\thanks{\footnotesize This work was
supported by Korean Ministry of Education, Sciences and Technology through NRF grant No. 2010-0017532 (H.K.) and 2011-0009671 (K.Y.).}}

\author{Hyeonbae Kang\thanks{\footnotesize Department of Mathematics, Inha University, Incheon 402-751, Republic of Korea (hbkang@inha.ac.kr)} \and Mikyoung Lim\thanks{Department of Mathematics, Korea Advanced Institute of Science and Technology, Yuseong-gu, Daejeon 305-701, Republic of Korea (mklim@kaist.ac.kr)} \and KiHyun Yun\thanks{\footnotesize Department of Mathematics,
Hankuk University of Foreign Studies,  Youngin-si, Gyeonggi-do 449-791, Republic of Korea (gundam@hufs.ac.kr)}}


\maketitle

\begin{abstract}
When two perfectly conducting inclusions are located closely to each other, the electric field concentrates in a narrow region in between two inclusions, and becomes arbitrarily large as the distance between two inclusions tends to zero. The purpose of this paper is to derive an asymptotic formula of the concentration which completely characterizes the singular behavior of the electric field, when inclusions are balls of the same radii in three dimensions.
\end{abstract}

\noindent {\footnotesize {\bf Mathematics subject classification
(MSC2000): 35J25, 73C40} }


\section{Introduction and statement of results}

Let $D_1$ and $D_2$ be bounded, simply connected and convex domains in $\Rbb^d$, $d=2,3$. Suppose that the conductivity of the inclusions is $\infty$, in other words, inclusions are perfect conductors. We consider the following conductivity problem:
\beq\label{eq_main}
\quad \left\{
\begin{array}{ll}
\ds\Delta u=0 \quad& \mbox{in } \mathbb{R}^d \setminus \overline{D_1 \cup D_2}, \\
\ds u=  C_j \ (\mbox{constant}) \quad&\mbox{on } \p D_j, j=1,2,\\
\ds u(\Bx) - H(\Bx)  =O(|\Bx|^{1-d}) \quad&\mbox{as } |\Bx| \to \infty,
\end{array}
\right.
\end{equation}
where $H$ is a given harmonic function in $\Rbb^d$ so that $-\nabla H$ is the background electric field in the absence of the inclusions. The constant value $C_j$ on $\p D_j$ is determined by the condition
\beq
\int_{\p D_j} \pd{u}{\nu^{(j)}} ~d\Gs=0 \quad \mbox{for }j=1,2.
\eeq
Here and throughout this paper $\nu^{(j)}$ is the outward unit normal to $\p D_j$.

The gradient of the solution $\nabla u$ represents the electric field (with the opposite sign) in the presence of inclusions and the stress field in two dimensional anti-plane elasticity, and it may become arbitrarily large as the distance between two inclusions tends to $0$. It has been proved that the generic rate of the gradient blow-up is $\Ge^{-1/2}$ in two dimensions \cite{AKLLL, AKLim, BLY, BC, keller,  Y, Y2} and $|\Ge \log \Ge|^{-1}$ in three dimensions \cite{BLY,
BLY2, lekner, LY}, where $\Ge$ is the distance between two inclusions. Occurrence of the gradient blow-up depends on the background potential (the harmonic function $H$
in \eqnref{eq_main}) and those background potentials which actually
make the gradient blow up are characterized in \cite{AKLLZ} when
$D_1$ and $D_2$ are disks.

The results mentioned above are estimates of the gradient of the solution from above and below, namely,
\beq\label{fracC_1}
\frac{C_1}{\psi(\Ge)} \le |\nabla u| \le \frac{C_2}{\psi(\Ge)} + C_3
\eeq
for some positive constants $C_1$, $C_2$ and $C_3$ where
\beq\label{psi}
\psi(\Ge) = \left\{
\begin{array}{ll}
\sqrt{\Ge} \quad &\mbox{if } d=2, \\
\ds \Ge \log \frac{1}{\Ge} \quad &\mbox{if } d=3.
\end{array}
\right.
\eeq
The constants $C_1$ and $C_2$ can possibly be $0$ depending on the background potential $H$.

The interest of this paper lies in the asymptotic behavior of $\nabla u$ as the distance between two inclusions tends to $0$. Since the singular behavior of $\nabla u$ occurs in the narrow region in between two inclusions, we are particularly interested in its behavior there.
In this regards, a complete characterization of the singular behavior of $\nabla u$ has been obtained when inclusions are disks \cite{KLY} and strictly convex domains in $\Rbb^2$ \cite{ACKLY}. Let $D_1$ and $D_2$ be disks in $\Rbb^2$ of radii $r_1$ and $r_2$, respectively, and let $R_j$ be the reflection with respect to $\p D_j$, $j=1,2$. Then the combined reflections $R_1R_2$ and $R_2R_1$ have unique fixed points, say $\Bf_1 \in D_1$ and $\Bf_2 \in D_2$.
Let
\beq\label{hdisks}
h(\Bx)=\frac{1}{2\pi}(\log |\Bx-\Bf_1|-\log |\Bx-\Bf_2|)
\eeq
(see section \ref{sec:sing} for a discussion on the function $h$). It has been proved that the solution $u$ to \eqnref{eq_main} can be expressed as
\beq\label{blowdisk}
u(\Bx)=\frac{4 \pi r_1r_2}{r_1+r_2}(\Bn\cdot\nabla H)(\Bc)h(\Bx)+g(\Bx), \quad \Bx\in\Rbb^2\setminus(D_1\cup D_2),
\eeq
where $\Bc$ is the middle point of the shortest line segment connecting $\p D_1$ and $\p D_2$, $\Bn$ is the unit vector in the direction of $\Bf_2-\Bf_1$, and $|\nabla g(\Bx)|$ is bounded independently of $\ep$ on any bounded subset of $\Rbb^2\setminus(D_1\cup D_2)$. So the singular behavior of $\nabla u$ is completely characterized by $\nabla h$. In particular, it can be shown using \eqnref{blowdisk} that the maximal concentration of $\nabla u$ occurs along the shortest line segment connecting $\p D_1$ and $\p D_2$, and on that segment
\beq\label{concendisk}
\nabla u \approx \frac{2\sqrt{2}}{\sqrt{\Ge}} \sqrt{\frac{r_1r_2}{r_1+r_2}} (\Bn\cdot\nabla H)(\Bc) \Bn.
\eeq

A complete characterization of the gradient blow-up like \eqnref{blowdisk}  has been obtained in  \cite{ACKLY} in the case when inclusions are strictly convex domains in $\Rbb^2$ by using disks osculating to convex domains. It is worth mentioning that the stress concentration factor for the $p$-Laplacian was derived in \cite{GN}.

The purpose of this paper is to derive an asymptotic formula of $\nabla u$ which characterizes its singular behavior when $D_1$ and $D_2$ are balls of the same radii in three dimensions.

In order to state the main result of this paper in a precise manner, let us fix notation.
Let $D_1$ and $D_2$ be balls of radius $r$  in three dimensions and $\Bc_1$ and $\Bc_2$ their centers. Let $\Bc$ be the middle point of $\Bc_1$ and $\Bc_2$, and $\Bn$ the unit vector in the direction of $\Bc_2-\Bc_1$, {\it i.e.},
$$
\Bc = \frac{\Bc_1+\Bc_2}{2}, \quad
\Bn= \frac{\Bc_2-\Bc_1}{|\Bc_2-\Bc_1|}.
$$
Let $R_j$, $j=1,2$, be the reflection with respect to $\p D_j$, \textit{i.e.},
$$
R_j(\mathbf{x})=\frac{r(\Bx-\Bc_j)}{|\Bx-\Bc_j|^2}+\Bc_j,
$$
and let, for $k=0,1, \ldots$,
\beq\label{repeat}
\begin{cases}
\ds\mathbf{p}_{2k}=(R_2R_1)^{k}\mathbf{c}_2, \\
\ds \mathbf{p}_{2k+1}= R_2(R_1R_2)^{k}\mathbf{c}_1.
\end{cases}
\eeq
We emphasize that $\Bp_n$ is contained in $D_2$ and monotonically converges to $\Bp$ as $n \to \infty$ where $\Bp$ is the fixed point of the combined reflection $R_2R_1$. Let
\beq
\mu_{n}=\frac{1}{|\mathbf{c}_1-\Bp_n|}, \quad n=1,2,\ldots,
\eeq
and
\beq
q_{0}=1 \quad \mbox{and} \quad q_{n}= \prod_{j=1}^n\mu_{j},\ n\geq 1 .
\eeq
Let $\Gr(\Bx)$ be the distance from $\Bx$ to the line connecting $\Bc_1$ and $\Bc_2$, {\it i.e.},
\beq\label{rhoax}
\Gr(\Bx)= |(\Bx-\Bc)-\la\Bx-\Bc,\mathbf{n}\ra \mathbf{n} |.
\eeq

The following is the main result of this paper.
\begin{thm}\label{main}
Suppose that the radius of the balls is much larger than the distance between them, {\it i.e.}, $\Ge \ll r$. The gradient $\nabla u$ of the solution to \eqref{eq_main} can be expressed as
\beq\label{asymp1}
\nabla u (\Bx) =  \frac {C_H^\Ge}{ \ds \left|\log \Ge \right| \left( \Ge + r\rho(\Bx)^2 \right)} (\mathbf{n}+ \eta(\Bx) )+ \nabla g (\Bx) \quad \mbox{if } \Gr(\Bx) \le \frac{r}{\ds \left|\log \Ge \right|^2}
\eeq
where
\beq\label{Cep}
C_H^\Ge = 2\sum_{n=0}^{\infty} q_n \left( H (\Bp_n)- H (-\Bp_n) \right),
\eeq
$|\nabla g|$ is bounded on any bounded region in $\Rbb^3 \setminus (D_1\cup D_2)$ regardless of $\Ge$, and
\beq
|\eta(\Bx)| \leq C |\log \Ge|^{-1}
\eeq
for some constant $C >0$ independent of $\Ge$.
\end{thm}

Some remarks on Theorem \ref{main} are in order. We first observe that the set $\Gr(\Bx) \le r|\log \Ge |^{-2}$ where \eqnref{asymp1} holds is a narrow region in between two spheres. The formula \eqnref{asymp1} shows that the major singular term of $\nabla u$ is
in the direction of $\Bn$, and that if $\Gr(\Bx)=\mbox{constant}$, then intensity of the field is constant. Note that the level set where $\Gr(\Bx)$ is constant is a cylinder around the line connecting centers of two spheres. So the intensity of the field decreases radially from the line connecting two centers of spheres.
The highest concentration of the field occurs when $\Gr(\Bx)=0$, in other words, when $\Bx$ is on the line segment connecting two closest points on the spheres, and on the segment,
\beq
\nabla u \approx \frac{C_H^\Ge}{\Ge |\log\Ge|} \Bn.
\eeq
Note that $C_H^\Ge$ depends on $\Ge$ since $\Bp_n$ and $q_n$ do. The following theorem reveals the limiting behavior of $C_H^\Ge$ as $\Ge \to 0$.
\begin{thm}\label{main2}
We have
\beq\label{0.3}
C_H^\Ge = C_H +O\left(\sqrt \epsilon |\log \epsilon|\right) \quad\mbox{as } \Ge \to 0
\eeq
where 
\beq
C_H = 2{\sum_{n=0} ^{\infty}} \frac 1 n \left( H\left(\frac r n \mathbf{n}+\mathbf{c}\right) -  H\left(-\frac r n \mathbf{n}+\mathbf{c}\right)\right).
\eeq
In particular, if $\Gr(\Bx)=0$, then
\beq\label{intensity}
\lim_{\Ge \to 0} \Ge |\log \Ge| |\nabla u(\Bx)| = |C_H|.
\eeq
\end{thm}

We emphasize that the occurrence of the gradient blow-up depends on the constant $C_H$: if $C_H \neq 0$, then it occurs. If $C_H=0$, then either $|\nabla u|$ is bounded or the blow-up rate is weaker than the generic rate $(\Ge|\log\Ge|)^{-1}$.  One can show for example that if the centers of the balls lie on the $x$-axis and their middle point is $(0,0,0)$, and if $H(x,y,z)=x^3-3xy^2$, then $C_H \neq 0$ and hence $|\nabla u|$ blows up as $\Ge \to 0$. It is interesting to observe that this is in contrast with two dimensional circular case. In view of \eqnref{concendisk}, the blow-up occurs only when $(\Bn\cdot\nabla H)(0,0) \neq 0$ (assuming $\Bc=(0,0)$). So, $\nabla u(x,y)$ blows up in two dimensions only when the background potential $H$ has the linear term $\Bn \cdot \Bx$.

The main ingredient in deriving \eqnref{asymp1} is the singular function $h$ which is the solution to
\beq\label{h:eqn}
\quad \left\{
\begin{array}{ll}
\ds\Delta h=0 \quad& \mbox{in } \mathbb{R}^d \setminus \overline{D_1 \cup D_2}, \\
\ds h= \mbox{constant}\quad& \mbox{on }\p D_j, \ j=1,2,\\
\nm
\ds \int_{\p D_j} \pd{h}{\nu^{(j)}} ~ds=(-1)^{j+1},\ \quad& j=1,2,\\
\ds h(\Bx)=O(|\Bx|^{1-d}) \quad&\mbox{as } |\Bx| \to \infty.
\end{array}
\right.
\end{equation}
Such a solution exists and is unique (see \cite{ACKLY, Y}). We emphasize that the constant values of $h$ on $\p D_1$ and on $\p D_2$ are different, and because of that the gradient of $h$ becomes arbitrarily large if the distance between $D_1$ and $D_2$ is small. This function characterizes the singular behavior of the solution to \eqnref{eq_main}. In fact, if we define the function $g$ by
\begin{equation}\label{old_result}
u(\Bx)  = \frac {u|_{\p D_2}- u|_{\p D_1} }{h|_{\p D_2}- h|_{\p D_1}}  h(\Bx) +  g(\Bx), \quad \Bx \in \Rbb^d \setminus (D_1 \cup D_2),
\end{equation}
then one can see that $g$ is harmonic in $\Rbb^d \setminus \overline{D_1 \cup D_2}$ and $g|_{\p D_1}= g|_{\p D_2}$, in other words, there is no potential difference of $g$ on $\p D_1$ and $\p D_2$. So it can be shown in the same way as in \cite{KLY} that $|\nabla g|$ is bounded on bounded subsets of $\Rbb^d \setminus (D_1 \cup D_2)$. It means that the singular behavior of $\nabla u$ is completely determined by $\frac {u|_{\p D_2}- u|_{\p D_1} }{h|_{\p D_2}- h|_{\p D_1}}  \nabla h(\Bx)$. Moreover, it is proved in \cite{Y, Y2} that
\beq\label{h:basic}
\ds u|_{\p D_1}- u|_{\p D_2} =  \int_{\p D_1}  H  \frac{\p h}{\p \nu^{(1)}} d\Gs + \int_{\p D_2}  H  \frac{\p h}{\p \nu^{(2)}} d\Gs,
\eeq
which means that the potential difference of $u$ is determined by the singular function $h$ (and the background potential $H$).

The function $h$ was first introduced in \cite{Y} and used in a crucial way to derive estimates for the gradient blow-up in \cite{LY, Y, Y2}. It is worth mentioning that $(\pd{h}{\nu^{(1)}}, \pd{h}{\nu^{(2)}})$ is an eigenvector corresponding to the eigenvalue $1/2$ of the Neumann-Poincar\'e operator associated with the interface problem \eqnref{eq_main} as shown in \cite{ACKLY, BT}.

If $D_1$ and $D_2$ are disks, then $h$ is given by \eqnref{hdisks}. In fact, $\p D_1$ and $\p D_2$ are the Apollonian circles of the fixed points $\Bf_1$ and $\Bf_2$, and hence $|\Bx-\Bf_1|/|\Bx-\Bf_2|$ is constant on $\p D_1$ and $\p D_2$. It is worth emphasizing that here the radii of disks may be different. If $D_1$ and $D_2$ are spheres, it is proved in \cite{LY} that $h$ is given by a weighted sum of the difference of the point charges: let $\GG(\Bx)= \frac{1}{4\pi} |\Bx|^{-1}$, the fundamental solution of the Laplacian in three dimensions. Then the singular function $h$ is given by
\begin{equation}\label{h:sameradii}
h(\Bx)=\frac{1}{\sum_{n=0}^{\infty} q_n} \sum_{n=0}^{\infty} q_n \left(\GG(\Bx - \Bp_n) - \GG(\Bx + \Bp_n) \right).
\end{equation}
This formula has been used in \cite{LY} to derive estimates like \eqnref{fracC_1}. We emphasize that in \cite{LY} an upper bound for $h$ is derived in a more general case when the radii of spheres are allowed to be different. In this paper we derive finer estimates of $h$ for the purpose of deriving \eqnref{asymp1}.

This paper is organized as follows. In section \ref{sec:sing}, we review the construction of the singular function in \cite{LY}. In section \ref{properties:pq}, we prove some technical lemmas which are required to estimate the singular function. In section \ref{sec:asymp}, we derive an asymptotic formula of the singular function. In the last section, we prove Theorem \ref{main} and Theorem \ref{main2}.

\section{Singular functions on spheres}\label{sec:sing}

Since the radius $r$ is much larger than $\Ge$, we may assume after scaling if necessary that $r=1$. We may also assume the centers are on the $x$-axis and $\Bc=(0,0,0)$ after rotation and shifting if necessary. We assume so in the sequel. It is also convenient to write $\Ge=2\Gd$ so that
$\Bc_1=(-1-\Gd, 0, 0)$ and $\Bc_2=(1+\Gd, 0, 0)$. Then, the function $\Gr$ defined in \eqnref{rhoax} becomes
\beq\label{rho2}
\rho(x,y,z)= \sqrt{y^2+z^2},
\eeq
and $\Bn = (1, 0, 0)$. Note that $\Bp_n$ defined by \eqnref{repeat} satisfies
\beq
\begin{cases}
\ds\mathbf{p}_{2k}=(R_2R_1)^{k}\mathbf{c}_2= - (R_1R_2)^{k}\mathbf{c}_1, \\
\ds \mathbf{p}_{2k+1}=-R_1(R_2R_1)^{k}\mathbf{c}_2 = R_2(R_1R_2)^{k}\mathbf{c}_1.
\end{cases}
\eeq

Define the function $h_1$ by
\beq
h_1(\Bx):=\sum_{k=0}^\infty\left(\frac{ q_{2k}}{|\Bx+\Bp_{2k}|}-\frac{ q_{2k+1}}{|\Bx - \Bp_{2k+1}|}\right) .
\eeq
Then $h_1$ is harmonic in $\Rbb^3 \setminus \overline{D_1 \cup D_2}$. Since the circle of Apollonius implies
\beq\label{Acircle}
|\By-\Bc_j||\Bx-R_j(\By)|= |\Bx-\By| \quad\mbox{for } |\By-\Bc_j|>1, \ \Bx\in\p D_j, \ j=1,2,
\eeq
we have
$$
\frac{ q_{2k+1}}{|\Bx - \Bp_{2k+1}|} = \frac{ q_{2k+1}}{|\Bc_1 - \Bp_{2k+1}|} \frac{1}{|\Bx- R_1 (\Bp_{2k+1})|} = \frac{ q_{2k+2}}{|\Bx + \Bp_{2k+2}|}
$$
if $\Bx \in \p D_1$, and
$$
\frac{ q_{2k}}{|\Bx + \Bp_{2k}|} = \frac{ q_{2k}}{|\Bc_2 + \Bp_{2k}|} \frac{1}{|\Bx- R_2 (-\Bp_{2k})|} = \frac{ q_{2k+1}}{|\Bx - \Bp_{2k+1}|}
$$
if $\Bx \in \p D_2$. So we have
\beq\label{h1const}
h_1|_{\p D_1}=1,\quad h_1|_{\p D_2}=0.
\eeq
Moreover, since
$$
\frac{1}{4\pi} \int_{\p D_j} \frac{\p}{\p \nu_\Bx} \frac{1}{|\Bx-\By|} d\Gs(\Bx)=
\begin{cases}
-1 \quad &\mbox{if } \By \in D_j, \\
0 \quad &\mbox{if } \By \notin \overline{D_j},
\end{cases}
$$
we have
\beq\label{h1int}
\frac{1}{4\pi} \int_{\p D_1} \frac{\p h_1}{\p \nu} d\Gs = -\sum_{k=0}^\infty q_{2k}, \quad \frac{1}{4\pi}\int_{\p D_2} \frac{\p h_1}{\p \nu} d\Gs = \sum_{k=0}^\infty q_{2k+1}.
\eeq

Define $h_2$ by
\beq
h_2(\Bx):=\sum_{k=0}^\infty\left(\frac{ q_{2k}}{|\Bx-\Bp_{2k}|}-\frac{ q_{2k+1}}{|\Bx + \Bp_{2k+1}|}\right) .
\eeq
Then $h_2$ is harmonic in $\Rbb^3 \setminus \overline{D_1 \cup D_2}$, and one can show similarly that
\beq\label{h2const}
h_2|_{\p D_1}=0,\quad h_2|_{\p D_2}=1,
\eeq
and
\beq\label{h2int}
\frac{1}{4\pi} \int_{\p D_1} \frac{\p h_2}{\p \nu} d\Gs = \sum_{k=0}^\infty q_{2k+1}, \quad \frac{1}{4\pi} \int_{\p D_2} \frac{\p h_2}{\p \nu} d\Gs = -\sum_{k=0}^\infty q_{2k}.
\eeq

It then follows from \eqnref{h1const}, \eqnref{h1int}, \eqnref{h2const}, and \eqnref{h2int} that
the solution to \eqnref{h:eqn} is given by
\beq\label{h:sameradii2}
h(\Bx):= -\frac{1}{4\pi \sum_{n=0} ^{\infty} q_n} \Bigr({h_1(\Bx)-h_2(\Bx)}\Bigr) = \frac{1}{4\pi \sum_{n=0} ^{\infty} q_n}
\sum_{n=0}^{\infty} q_n \left(\frac 1 {|\Bx - \Bp_n|} - \frac 1 {|\Bx + \Bp_n|} \right).
\eeq
Thus we have \eqnref{h:sameradii}. We also have
\beq\label{hconst}
h|_{\p D_2}- h|_{\p D_1} = \frac{2}{4 \pi\sum_{n=0}^{\infty} q_n}.
\eeq

In the next section we derive fine properties of the sequences $\Bp_n$ and $q_n$, which are used in deriving an asymptotic formula for $h$.

\section{Properties of the sequences $\Bp_n$ and $q_n$}\label{properties:pq}

Let $\Bp=(p,0,0)$ be the fixed point of the combined reflection $R_2R_1$ as before. Then one can easily see that
$p$ satisfies
$$
p = -\frac 1 {1+\Gd +p} +1 +\Gd,
$$
so that
\beq\label{sqrt2d}
p  =  \sqrt{2\Gd}  + O (\Gd) \quad\mbox{as }\Gd \to 0.
\eeq
Let $\Bp_n=(p_n, 0, 0)$. Then, $p_0=1+\Gd$ and $p_n$ satisfies the recursive relations
\beq\label{p_n+1}
p_{n+1} = -\frac 1 {1+\Gd +p_n} +1 +\Gd, \quad n=0, 1, \ldots.
\eeq
One can further see that
\beq\label{pnp}
p_n = p \left(\frac {2} {A^{n+1}-1} + 1\right)=p \left(\frac {A^{n+1}+1} {A^{n+1}-1} \right),
\eeq
where
\beq\label{def:A}
A:= \frac {1 + \Gd + p }{1+ \Gd - p } .
\eeq
Note that
\beq\label{def:A2}
A = 1 + 2p + O(\Gd) = 1+ 2\sqrt{2\Gd} + O(\Gd).
\eeq
In particular, the sequence $p_n$ is decreasing and converges to $p$ as $n \to \infty$.

Since
\beq
\mu_{n}= \frac{1}{|\mathbf{c}_1-\Bp_n|} = \frac 1 {1+\Gd+p_n}=(1+\Gd -p_{n+1}),
\eeq
we have
\beq\label{q_n+1}
q_{n+1} = \mu_{n} q_n = \frac 1 {1+\Gd+p_n} q_n =(1+\Gd -p_{n+1}) q_n.
\eeq

For a given $\Gd>0$, let $N_0=N_0(\Gd)$, $N=N(\Gd)$ and $N_1=N_1(\Gd)$ be as follows:
\beq\label{notation}
N_0(\Gd) = \left[{|\log \Gd|}\right], \quad N(\Gd)=\left[\frac 1 {\sqrt {\Gd}}\right], \quad N_1(\Gd)=\left[\frac {1}  {\Gd |\log \Gd|}\right].
\eeq
Here $[\cdot]$ is the Gaussian bracket. We use this notation for the rest of this paper. Since $\Gd$ is sufficiently small, we have
$$N_0(\Gd)\ll N(\Gd) \ll N_1(\Gd).$$

The following lemma was obtained in \cite{LY}.
\begin{lem}\label{lem:pq1overn}
There is a constant $C$ independent of $\Gd$ such that
\beq\label{pq1overn}
\left| p_n - \frac 1 {n+1} \right| + \left| q_n - \frac 1 {n+1} \right| \le C \sqrt \Gd
\eeq
and
\beq
|p_n - p_{n+1}| < \frac C {n^2}
\eeq
for $n \leq N(\Gd)$.
\end{lem}

We prove the following lemma.
\begin{lem} \label{elem-propty}
Let $N=N(\Gd)$ and $N_1=N_1(\Gd)$ as before.
\begin{itemize}
\item[\rm(i)] There is a positive $C$ independent of $\Gd$ such that
\beq
\left|\sum_{n=0}^{\infty}  q_n  - \sum_{n=1} ^{N}  \frac 1 n   \right| \leq C \quad \mbox{and} \quad
\sum_{n=N}^\infty q_n\leq C.
\eeq
\item[\rm(ii)] $p_n -p \geq 2  \sqrt {\Gd} A^{-n}$ for all $n$.
\item[\rm(iii)] There is a constant $C$ such that
$$
p_n-p\geq\frac{C}{n}
$$
for all $n \leq N$.
\item[\rm(iv)] $0< p_{N_1} - p \leq e^{-1 /(\sqrt \Gd | \log \Gd|)}$.
\end{itemize}
\end{lem}
\pf
Since $p_n$ decays to $p$, we have from \eqnref{q_n+1}
\beq\label{qnqm}
q_n \leq (1+\Gd-p)^{n-m} q_m \quad \mbox{for all }n\geq m \geq 1 .
\eeq
So, it follows from \eqnref{pq1overn} that
$$
\sum_{n=N}^\infty q_n \le \sum_{n=N}^{\infty} q_N (1+\Gd - p)^{n- N} \le \left(C\sqrt{\Gd} + \frac{1}{N+1}\right) \sum_{n=N}^{\infty} (1+\Gd - p)^{n- N} \le C,
$$
and
$$
\left|\sum_{n=0} ^{\infty}  q_n  - \sum_{n=1} ^{N}  \frac 1 n  \right| \leq
C ~N { \sqrt \Gd}+ \sum_{n=N}^{\infty} q_n  \leq C.
$$
This proves (i).

We have from \eqnref{pnp} that for each $n\in\mathbb{N}$,
\beq\label{pn-p}
p_n-p=\frac {2p}{A^{n+1}-1 } .
\eeq
So, (ii) follows from \eqnref{sqrt2d}.

Now, suppose that $n\leq N$. Since $A \le 1+3p$, using the inequality
$$
(1+s)^{n} \leq  1 + ns + \frac 1 2 n^2 s^2 (1+s)^n
$$
which holds for all $s>0$, we obtain
\beq\label{Anle}
A^{n}\le (1+3p)^{n} \leq 1 +  3n p + \frac{9}{2} n^2 p^2 (1+3p)^{n}.
\eeq
Since $np \le Np \le 2$ and $(1+t)^{1/t} $ increases to $e$ as $t \rightarrow 0+$, we have
$$
(1+3p)^{n} \leq \left( (1 + 3p)^{\frac 1 {3p}}\right)^{3np} \le e^{6} ,
$$
and hence, from the second inequality in \eqnref{Anle}
$$
A^n \le 1 + C np
$$
for some constant $C$ independent of $n \le N$ and $\Gd$.
We then infer from \eqnref{pnp} that
$$
p_n - p = \frac{2p}{ A^{n+1}-1}  \geq \frac 1 {Cn} ,\quad n\leq N(\Gd).
$$

Now, if $n=N_1$, then we have
$$
\log(A^n) = n \log A \ge \frac{n(A-1)}{2} \ge \frac{1}{\sqrt\Gd|\log\Gd|},
$$
and hence
$$
A^n\geq e^{\frac{1}{\sqrt\Gd|\log\Gd|}}.
$$
Now (iv) follows from \eqnref{pn-p}. This completes the proof.
\qed

Lemma \ref{elem-propty} (i) yields
\beq\label{q:sum}
\ds\sum_{n=0} ^{\infty}  q_n =\frac{1}{2}{|\log\Gd|}+O(1).
\eeq

The following lemma provides the finer properties of $p_n$ and $q_n $ that are crucial in proving the main result of this paper.

\begin{lem}\label{key_lemma}
\begin{itemize}
\item[\rm(i)] If $ N_0(\Gd) \leq n \leq N_1(\Gd)$, then
\beq
\frac {q_n}{p_n - p_{n+1}} = \frac {1+ O(|\log \Gd|^{-1})}{\sqrt {p_n ^2 - p^2}}\quad \mbox{as }\Gd \rightarrow 0,
\eeq
where $O(|\log \Gd|^{-1})$ is independent of $n$.

\item[\rm(ii)] There are constants $C_1$ and $C_2$ such that
\beq
q_{n} \leq C_1   (1- p + \Gd)^{n- N_1} e^{- \frac {C_2}{\sqrt {\Gd}|\log \Gd|}}
\eeq
for all $n\geq N_1=N_1(\Gd)$.

\end{itemize}
\end{lem}
\pf
If $n>m$, then we have from \eqnref{q_n+1}
$$
\log q_n = - \sum_{j=m}^{n-1} \log (1+\Gd + p_{j}) + \log q_m.
$$
Using the inequality $|\log(1+t)-t| \le C t^2$, we obtain
\begin{align*}
\log q_n = - \sum_{j=m}^{n-1} p_j - \Gd (n-m) + \log q_m +E_1,
\end{align*}
where the error term $E_{1}$ satisfies
\beq
|E_{1}| \leq  C_1 \sum_{j=m}^{n-1} (\Gd + p_j)^2 \le C_2 \sum_{j=m}^{n-1} p_j ^2.
\eeq
The last inequality above holds since $ \Gd \ll p < p_j $. Here and in the rest of this proof, $E_j$'s denote errors to be estimated. We then have from \eqnref{pnp} that
\begin{align*}
\log \frac{q_n}{q_m} = -(\Gd+ p) (n-m)- 2p \sum_{j=m+1}^{n} \frac {1}{A^{j}-1}  + E_{1} .
\end{align*}
Since $\frac {1}{A^j-1}=\frac{A^{-j}}{1-A^{-j}}$ is decreasing in $j$, we have
\begin{align*}
& \left| \sum_{j=m+1}^{n} \frac {1}{A^j-1} + \frac{1}{\log A} \log \left ( \frac{1-A^{-m-1}}{1-A^{-n-1}} \right) \right| \\
& = \left| \sum_{j=m+1}^{n} \frac{A^{-j}}{1-A^{-j}} - \int_{m+1}^{n+1} \frac{A^{-x}}{1-A^{-x}}dx \right| \le
\frac{A^{-m-1}}{1- A^{-m-1}}.
\end{align*}
So we have
\begin{align*}
\log \frac{q_n}{q_m} = -(\Gd+ p) (n-m) + \frac{2p}{\log A} \log \left ( \frac{1-A^{-m-1}}{1-A^{-n-1}} \right) + E_{2},
\end{align*}
where the new error term $E_2$ satisfies
\beq\label{E_2}
|E_2| \le C \left( \sum_{j=m}^{n-1} p_j ^2 + \frac{pA^{-m-1}}{1- A^{-m-1}} \right).
\eeq

One can see from \eqnref{def:A2} that
\beq
\frac{2p}{\log A}=1+ E_3,
\eeq
where
\beq\label{E_3}
|E_3| \le C \sqrt{\Gd}.
\eeq
So, we have
\begin{align*}
\log \frac{q_n}{q_m} = -(\Gd+ p) (n-m) + (1+E_3) \log \left ( \frac{1-A^{-m-1}}{1-A^{-n-1}} \right) + E_{2} ,
\end{align*}
which in turn implies
\beq\label{q_n/q_m}
q_n =  q_m e^{-pn} \left ( \frac{1-A^{-m-1}}{1-A^{-n-1}} \right) e^{E_4},
\eeq
where
\beq\label{E_4}
E_4:= pm - \Gd(n-m) + E_2 + E_3 \log \left ( \frac{1-A^{-m-1}}{1-A^{-n-1}} \right).
\eeq

Note that
$$
p_n - p_{n+1} = (p_{n+1}-p) \frac {A-1}{1-A^{-n-1}},
$$
so that
\beq\label{eq:q_p_m-p_m}
\frac{q_n}{p_n - p_{n+1}} =
\frac {q_m e^{-p n}}{p_{n+1}-p} \frac{1-A^{-m-1}}{A-1} e^{E_4} .
\eeq
Since $p_n/p = (A^{n+1}+1)/(A^{n+1}-1)$, we have
$$
(n+1) \log A = \log \left( \frac { p_n+p } { p_n-p }\right) ,
$$
and, since $\log A= 2p + O(\Gd)$, it follows that
$$
p n = \frac{1}{2} \log \left( \frac { p_n+p } { p_n-p }\right) + E_5,
$$
where
\beq\label{E_5}
|E_5| \le C (n \Gd + \sqrt{\Gd}).
\eeq
We then obtain from \eqnref{eq:q_p_m-p_m}
\begin{align}
\frac{q_n}{p_n - p_{n+1}} &=  \sqrt{\frac{p_n-p}{p_n+p}}\frac {1}{p_{n+1}-p}
q_m \frac{1-A^{-m-1}}{A-1} e^{E_4 - E_5} \nonumber \\
&= \frac{1}{\sqrt{p_n^2-p^2}}\frac {p_n-p}{p_{n+1}-p}
q_m \frac{1-A^{-m-1}}{A-1} e^{E_4 - E_5}
.  \label{eq:q_p_m-p_m2}
\end{align}

Suppose now that $m= N_0-1$ and $m<n\leq  N_1$. Then we have $E_5=O(|\log \Gd|^{-1})$. We will show that
\begin{align}
& \frac {p_n-p}{p_{n+1}-p} = 1+ O(|\log \Gd|^{-1}), \label{estimate1} \\
& q_m \frac{1-A^{-m-1}}{A-1} = 1+ O(|\log \Gd|^{-1}), \label{estimate2} \\
& E_4 = O(|\log \Gd|^{-1}).  \label{estimate3}
\end{align}
Once we have these estimates, then (i) follows from \eqnref{eq:q_p_m-p_m2}.

To prove \eqnref{estimate1}, we first observe that
\begin{align*}
\frac {p_n - p }{p_{n+1} - p} = \frac {A^{n+2}-1}{A^{n+1}-1} = A \left( 1 +  \frac {1}{A +A^2 +\cdots + A^{n+1}  }\right).
\end{align*}
Since $A>1$, $n \ge |\log\Gd|$ and $A=1+O(\sqrt{\Gd})$, we have
$$
\frac {p_n - p }{p_{n+1} - p} = (1+O(\sqrt{\Gd}))(1+ O(|\log \Gd|^{-1})) = 1+ O(|\log \Gd|^{-1}).
$$

To prove \eqnref{estimate2}, we use inequalities
$$
(m+1)s - \frac 1 2 m(m+1) s^2 \leq 1-(1-s)^{m+1} \leq  (m+1)s
$$
which hold for all $s \in [0,1]$. Since $A^{-1}=1-2p+ O(\Gd)$, we have
$$
(m+1)(2p+O(\Gd)) - \frac 1 2 m(m+1)(2p+O(\Gd))^2 \le 1-A^{-m-1} \le (m+1)(2p+O(\Gd)).
$$
Since $m=O(|\log\Gd|)$ and $p=O(\sqrt{\Gd})$, we have
\beq\label{1-Am}
1-A^{-m-1}=2(m+1)p+O(\Gd |\log\Gd|^{2}).
\eeq
Note that
$$
\frac{1}{A-1} =  \frac{1}{2p + O(\Gd)} = \frac{1}{2p} + O(1).
$$
Since $q_m = \frac{1}{m+1} + O(\sqrt{\Gd})$ by Lemma \ref{lem:pq1overn} and $m=N_0-1$,  we infer that
\begin{align*}
q_m \frac{1-A^{-m-1}}{A-1} & = \left( \frac{1}{m+1} + O(\sqrt{\Gd}) \right) \left( 2(m+1)p+O(\Gd |\log\Gd|^{2}) \right) \left( \frac{1}{2p} + O(1) \right)\\
& = 1+ O(|\log \Gd|^{-1}).
\end{align*}
So, \eqnref{estimate2} is proved.

To prove \eqnref{estimate3}, we first estimate $E_2$. We have from \eqnref{pq1overn} that
\begin{align*}
\sum_{j=m+1}^{n} p_j^2 &\leq C\sum_{j=N_0}^N p_j^2 + \sum_{j=N+1}^{N_1} p_j^2\\
&\leq C \sum_{j=N_0}^N \frac{1}{j^2} + \sum_N^{N_1}p_N^2\\
&\leq C \left( \frac{1}{N_0} + p_N^2 N_1 \right ) = O(|\log \Gd|^{-1}).
\end{align*}
On the other hand, it follows from \eqnref{1-Am} that
$$
\frac{pA^{-m-1}}{1- A^{-m-1}} = \frac{p + O(\Gd |\log\Gd|)}{2(m+1)p+O(\Gd |\log\Gd|^{2})} = \frac{1}{2(m+1)} \left(1 + O(\sqrt{\Gd} |\log\Gd|) \right) = O(|\log \Gd|^{-1}).
$$
So we infer from \eqnref{E_2} that
\beq\label{E_22}
E_2= O(|\log \Gd|^{-1}).
\eeq

Since $p=O(\sqrt{\Gd})$, we obtain from \eqnref{1-Am} that
$$
1\geq \frac {1-A^{-m-1}}{1-A^{-n-1}}\geq 1-A^{-m-1} \geq  C\sqrt{\Gd} |\log \Gd|.
$$
We then infer from \eqnref{E_3} that
\beq\label{E_32}
\left| E_3 \log \left ( \frac{1-A^{-m-1}}{1-A^{-n-1}} \right) \right| \le C \sqrt{\Gd} |\log \Gd|.
\eeq
Thus we have from \eqnref{E_4}, \eqnref{E_22} and \eqnref{E_32} that
\begin{align*}
|E_4| &\le |pm| + \Gd (n-m) + |E_2| + \left| E_3 \log \left ( \frac{1-A^{-m-1}}{1-A^{-n-1}} \right) \right| \\
& \le C \left( \sqrt{\Gd} |\log\Gd| + |\log\Gd|^{-1} \right) \le C |\log\Gd|^{-1},
\end{align*}
so \eqnref{estimate3} is proved.

We have from  \eqref{q_n/q_m} that
$$
q_{N_1} =  q_m e^{-pN_1} \left ( \frac{1-A^{-m-1}}{1-A^{-N_1-1}} \right) e^{E_4}.
$$
So, it follows from \eqnref{estimate3} that
$$
q_{N_1} \leq C_1  e^{- \frac {C_2}{\sqrt {\Gd}|\log \Gd|}}
$$
for some constants $C_1$ and $C_2$. Now,
(ii) follows from \eqnref{qnqm}. This completes the proof.
\qed

\section{Asymptotic behavior of the singular function}\label{sec:asymp}

Let
\beq
R_\Gd:= \left\{ ~ \Bx \in \Rbb^3 \setminus (D_1 \cup D_2)~ \big| ~ \Gr(\Bx) \le |\log \Gd|^{-2} ~ \right\} ,
\eeq
where $\Gr$ is given by \eqnref{rho2}. Note that $R_\Gd$ is a narrow region in between $D_1$ and $D_2$.
Let
\beq
v(\Bx):= \Bigr(4\pi \sum_{n=0}^{\infty} q_n \Bigr) h(\Bx)= \sum_{n=0}^{\infty} q_n \left(\frac 1 {|\Bx - \Bp_n|} - \frac 1 {|\Bx + \Bp_n|} \right).
\eeq
In this section we investigate the asymptotic behavior of $\nabla v(\Bx)$ in the region $R_\Gd$.
We obtain the following proposition.
\begin{prop}\label{main_prop}
For $\Bx=(x,y,z)\in R_\Gd$, we have
\beq\label{asymp_sing}
\nabla v(\Bx) = \frac{2}{2\Gd + \Gr(\Bx)^2}\left((1,0,0) + O(|\log \Gd|^{-1}) \right).
\eeq
\end{prop}

It turns out that $|\p_y v (\Bx)|$ and $|\p_z v(\Bx)|$ can be estimated without much difficulty. In fact, we obtain the following lemma whose proof is given in Subsection \ref{pf:prop3}.
\begin{lem}\label{prop3}
For $\Bx=(x,y,z)\in R_\Gd$, we have
\beq\label{h_r}
|\p_y v (\Bx)|+|\p_z v(\Bx)| \leq  \frac{C}{\sqrt \Gd + \Gr(\Bx)} \left( 1+\log\Bigr(1+\frac{\Gr(\Bx)^2}{\Gd}\Bigr) \right)
\eeq
for some constant $C$ independent of $\Gd$.
\end{lem}

Estimates of $\p_x v (\Bx)$, especially those terms for $N_0 \le n \le N_1$, are quite involved. Based on Lemma \ref{key_lemma} (i) we compare $v$ which is given as an infinite series with the integral defined by
\beq\label{def:h_0}
v_0(\Bx):= \int_{p} ^1 \left(\frac  1  {|\Bx - (t,0,0)|} - \frac  1  {|\Bx + (t,0,0)|}  \right)~ \frac 1 {\sqrt {t^2  - p^2 } } ~d t,
\eeq
where $(p,0,0)$ be the fixed point of the combined reflections $R_2R_1$. We obtain the following lemmas whose proofs are given in Subsection \ref{pf:prop4} and \ref{pf:prop2}, respectively.
\begin{lem}\label{prop4}
For $\Bx=(x,y,z)\in R_\Gd$, we have
\beq\label{prop4_formul}
\partial_x v_0(\Bx)= \frac 2 {2\Gd + \Gr(\Bx)^2}\left(1 + O(|\log \Gd|^{-1})  \right)  .
\eeq
\end{lem}

\begin{lem}\label{prop2}
For $\Bx=(x,y,z)\in R_\Gd$,  we have
\beq\label{d_xh}
\p_x v (\Bx) = \p_x v_0(\Bx)\left(1+O(|\log\Gd|^{-1})\right).
\eeq
\end{lem}

Proposition \ref{main_prop} is an immediate consequence of above lemmas.

\subsection{Proof of Lemma \ref{prop3}}\label{pf:prop3}

We first observe that if $\Bx = (x,y,z) \in R_\Gd$, then $|x|\leq 1+\Gd-\sqrt{1-y^2-z^2}$ and $\Gr \leq |\log \Gd|^{-2}$, and hence
\beq\label{xinRep}
|x|\leq \Gd +\Gr(\Bx)^2.
\eeq

Using notation $\Gr=\Gr(\Bx)$, $v$ can expressed as
\beq
v(\Bx) = \sum_{n=0}^{\infty} q_n \left(\frac{1}{\sqrt{(x-p_n)^2+ \Gr^2}} -  \frac{1}{\sqrt{(x+p_n)^2+ \Gr^2}} \right).
\eeq
So, it suffices to estimate $|\p_\Gr v|$. Note that
\begin{align*}
\p_\Gr \left(\frac{1}{\sqrt{(x-p_n)^2+ \Gr^2}} -  \frac{1}{\sqrt{(x+p_n)^2+ \Gr^2}} \right)
&= -\frac{\Gr}{[(x-p_n)^2+ \Gr^2]^{3/2}} + \frac{\Gr}{[(x+p_n)^2+ \Gr^2]^{3/2}} \\
&= 3 \Gr \int_{-x}^x \frac{t-p_n}{[(t-p_n)^2+ \Gr^2]^{5/2}} dt.
\end{align*}
Therefore, we have
$$
\left| \p_\Gr \left(\frac{1}{\sqrt{(x-p_n)^2+ \Gr^2}} -  \frac{1}{\sqrt{(x+p_n)^2+ \Gr^2}} \right) \right| \le
3 \Gr \int_{-x}^x \frac{1}{[(t-p_n)^2+ \Gr^2]^{2}} dt.
$$
By \eqnref{xinRep} we have
\beq\label{t-pn}
(t-p_n)^2+ \Gr^2 \ge C(\Gr^2+p_n^2)
\eeq
for some constant $C$. It then follows that
$$
\left| \p_\Gr \left(\frac{1}{\sqrt{(x-p_n)^2+ \Gr^2}} -  \frac{1}{\sqrt{(x+p_n)^2+ \Gr^2}} \right) \right| \le
C  \frac{\rho (\rho^2 + \Gd)}{\rho^4 + p_n^4} \le C  \frac{\rho}{\rho^2+ p_n^2}.
$$
So we have
$$
|\p_\Gr v(\Bx)| \le C \sum_{n=0}^\infty \frac{\rho q_n}{\rho^2+ p_n^2}.
$$

Let $N=N(\Gd)$. Using Lemma \ref{lem:pq1overn}, we have
\begin{align*}
\sum_{n=0}^{N-1} \frac{\rho q_n}{\rho^2+ p_n^2} & \leq  C \sum_{n=1}^{N}   \frac {\rho} {(1/n^2 + \rho^2)} ~ \frac 1 n \\
& \leq   C \sum_{n=1} ^{N} \frac {\rho n}  {1 + \rho^2n^2}  \leq  C \left(1+    \int_{1} ^{1/{\sqrt {\Gd}}} \frac {\rho s} {1 + \rho^2 s^2} ds\right)\\
& \leq  C \frac 1 {\rho + \sqrt \Gd} \left( 1+\log\Bigr(1+\frac{\rho^2}{\Gd}\Bigr) \right).
\end{align*}
If $n \ge N $, then $q_n =O(\sqrt{\Gd})$, and thus we have from \eqnref{qnqm}
\begin{align*}
\sum_{n=N}^{\infty} \frac{\rho q_n}{\rho^2+ p_n^2} \leq C \sum_{n=N}^{\infty} \frac {\rho} {\rho^2+\Gd} \sqrt {\Gd} (1+\Gd-p)^{n-N}
\leq \frac {C}{\rho+ \sqrt {\Gd}}.
\end{align*}
This completes the proof. \qed

\subsection{Proof of Lemma \ref{prop4}}\label{pf:prop4}

Let
\begin{align*}
\p_x v_0 (\Bx) =& \int_{p} ^{|\log \Gd|^{-1}} + \int_{|\log \Gd|^{-1}}^{1} \p_{x} \left(\frac  1  {|\Bx - (t,0,0)|} - \frac  1  {|\Bx + (t,0,0)|}  \right)~ \frac 1 {\sqrt {t^2  - p^2 } } ~d t \\
:= & \, I + II.
\end{align*}

If $|\log \Gd|^{-1} \le t \le 1$, then $|\Bx \pm (t,0,0)| \ge C t$ for some constant $C$ and for all $\Bx \in \Rbb^3 \setminus (D_1 \cup D_2)$. Since $p = O(\sqrt{\Gd})$, we also have $\sqrt {t^2  - p^2} \ge C t$. Thus we have
\beq\label{cbig}
|II| \leq  C \int_{|\log \Gd|^{-1}}^{1}~\frac 1 {t^3} ~d t  \leq  C |\log \Gd|^2 .
\eeq

Suppose now that $p\leq  t \leq |\log \Gd|^{-1}$. Using \eqnref{xinRep} and the fact that $p= O(\sqrt\Gd)$ again, we have for all $\Bx \in R_{\Gd}$
\begin{align}
&|tx|\leq t (\Gr^2 + \Gd)\leq \frac C {|\log \Gd|} (t^2+ \Gr^2) , \label{eq_lem1} \\
&|x|^2\leq (\Gr^2 + \Gd)^2\leq \frac C {|\log \Gd|} (t^2+ \Gr^2) \label{eq_lem2}
\end{align}
for some constant $C$ independent of $\Gd$. Thus, we have
\begin{align}
\frac  1  {|\Bx \pm (t,0,0)|^3} =&  \frac  1  {{\left((x\pm t)^2 + \Gr^2 \right)}^{3/2}}=  \frac  1  {{\left( t^2+ \Gr^2 \pm 2 xt + x^2 \right)}^{3/2}} \nonumber \\
=&  \frac{1}{( t^2+ \Gr^2 )^{3/2}} \left( 1+ O \left(|\log \Gd|^{-1} \right)  \right). \label{lemA1}
\end{align}

 From the mean value property, we have
\begin{align}
& \left| \frac  {-1}  {|\Bx - (t,0,0)|^3} + \frac  {1}  {|\Bx + (t,0,0)|^3}\right|\notag\\
&= \left|\frac  {1}  {\bigr((   t^2 + x^2 +\Gr^2) - 2xt \bigr)^{3/2}} - \frac  {1}  {\bigr(( t^2 + x^2 +\Gr^2 ) + 2xt \bigr)^{3/2}}\right|\notag\\
&\leq  \frac  {6|xt|}  {\left| (t^2 + x^2 +\Gr^2 ) - |2xt| \right|^{5/2}}.\label{mvpxmc}
\end{align}
It then follows from \eqref{eq_lem1} and \eqref{eq_lem2} that
\beq\label{lemA2}
\left|x \left(\frac  {-1}  {|\Bx - (t,0,0)|^3} + \frac  {1}  {|\Bx + (t,0,0)|^3}\right)\right| \leq C |\log \Gd|^{-1}
  \frac t {(t^2 + \Gr(\Bx)^2 )^{3/2}}
\eeq
for some constant $C$ independent of $\Gd$.

Since
\begin{align} &\p_{x} \left(\frac  1  {|\Bx - (t,0,0)|} - \frac  1  {|\Bx + (t,0,0)|}  \right)\nonumber\\\nonumber
& = -\frac  {x-t}  {|\Bx - (t,0,0)|^3} + \frac  {x+t}  {|\Bx + (t,0,0)|^3}  \\ & = t\left( \frac  1  {|\Bx - (t,0,0)|^3} + \frac  1  {|\Bx + (t,0,0)|^3} \right)+ x \left(\frac  {-1}  {|\Bx - (t,0,0)|^3} + \frac  {1}  {|\Bx + (t,0,0)|^3}\right), \label{pintegrand}
\end{align}
we obtain from \eqnref{lemA1} and \eqnref{lemA2}
$$
\p_{x} \left(\frac  1  {|\Bx - (t,0,0)|} - \frac  1  {|\Bx + (t,0,0)|}  \right) =
\frac{t}{( t^2+ \Gr^2 )^{3/2}} \left( 2+ O \left(|\log \Gd|^{-1} \right)  \right).
$$
It then follows that
$$
I= \left(2 + O\left(|\log \Gd|^{-1} \right)  \right)  \int_{p}^{|\log \Gd|^{-1}} \frac  t  {{\left( t^2+ \Gr^2 \right)}^{3/2}}
~ \frac 1 {\sqrt {t^2  - p^2 } } ~d t.
$$
Using the substitution $t=\sqrt{t^2-p^2}$, one can see that
$$
\int_{p}^{|\log \Gd|^{-1}} \frac  t  {{\left( t^2+ \Gr^2 \right)}^{3/2}}
~ \frac 1 {\sqrt {t^2  - p^2 } } ~d t
= \frac 1 {p^2 +\Gr^2}
\left( \frac{|\log \Gd|^{-2} - p^2}{|\log \Gd|^{-2}+\Gr^2} \right)^{1/2}.
$$
Since $p=\sqrt{2\Gd} + O(\Gd)$ and $\Gr \le |\log\Gd|^{-2}$, we have
$$
\frac 1 {p^2 +\Gr^2}
\left( \frac{|\log \Gd|^{-2} - p^2}{|\log \Gd|^{-2}+\Gr^2} \right)^{1/2} = \frac 1 {2\Gd +\Gr^2}
\left( 1+ O \left(|\log \Gd|^{-1} \right)  \right),
$$
and hence
$$
I = \frac 2 {2\Gd +\Gr^2}
\left( 1+ O \left(|\log \Gd|^{-1} \right)  \right),
$$
which together with \eqnref{cbig} yields
$$
\p_x v_0(\Bx) = \frac 2 {2\Gd +\Gr^2}
\left( 1+ O \left(|\log \Gd|^{-1} \right)  \right) + O(|\log \Gd|^2).
$$
Since $\Gr \le |\log\Gd|^{-2}$, the above formula can be written as
$$
\p_x v_0(\Bx) = \frac 2 {2\Gd +\Gr^2}
\left( 1+ O \left(|\log \Gd|^{-1} \right)  \right) .
$$
This completes the proof. \qed

\subsection{Proof of Lemma \ref{prop2}}\label{pf:prop2}

Let $N_0= [|\log\Gd|]$ and $N_1= [\frac{1}{\Gd|\log\Gd|}]$ as before. Let
\begin{align*}
\p_x v (\Bx) \notag = & \sum_{n=0}^{N_0-1}+\sum_{n=N_0}^{N_1-1} +\sum_{n=N_1}^{\infty} \p_x\left(  \frac 1 {|\Bx - \Bp_n|}  - \frac 1 {|\Bx + \Bp_n|} \right)q_n\\
:= & \, S_1(\Bx)+S_2(\Bx)+S_3(\Bx),
\end{align*}
and
\begin{align*}
\p_x v_0(\Bx) =& \int_{p_{N_0}}^{1} + \int_{p_{N_1 }} ^{p_{N_0}}+ \int_{p} ^{p_{N_1  }} \p_x\left(  \frac 1 {|\Bx - (t,0, 0)|}  - \frac 1 {|\Bx + (t,0, 0)|}  \right) \frac {1}{\sqrt {t^2 - p^2 }} dt \nonumber\\
  := &\, I_1+I_2+I_3. 
\end{align*}

We first estimate $S_1$, $S_3$, $I_1$, and $I_3$. There is a constant $C>0$ independent of $n$ such that $|\Bx \pm \Bp_n| \ge C p_n$ for all $\Bx \in R_\Gd$. So we have from \eqnref{pq1overn} that
\begin{align*}
 |S_1(\Bx)|&\leq \sum_{n=0}^{N_0-1}  \left|\nabla \left(\frac  1  {|\Bx - \Bp_n|} - \frac  1  {|\Bx + \Bp_n|}  \right) q_n \right| \\
 & \leq C \sum_{n=1}^{N_0-1}  n^2 ~\frac 1 n \leq C |\log \Gd|^2.
\end{align*}
We also have from Lemma \ref{key_lemma} (ii) that
\begin{align*}
 |S_3(\Bx)|&\leq \sum_{n=N_1}^{\infty}  \left| \nabla \left(\frac  1  {|\Bx - \Bp_n|} - \frac  1  {|\Bx + \Bp_n|}  \right) q_n \right|
 \le \sum_{n=N_1}^{\infty} \frac{1}{p^2} q_n  \\
 & \leq C \sum_{n=N_1}^{\infty}  \frac 1 {\Gd} ~ (1- p + \Gd)^{n- N_1} e^{- C_2 \frac {1}{\sqrt {\Gd}|\log \Gd|}} \leq C.
\end{align*}
Similarly, we have
\begin{align*}
|I_1(\Bx)| &\leq \int_{p_{N_0}} ^{1} \left| \nabla  \left(\frac  1  {|\Bx - (t,0,0)|} - \frac  1  {|\Bx + (t,0,0)|}  \right) ~\frac 1 {\sqrt {t^2- p^2}} \right| ~dt  \\
& \leq C \int_{|\log \Gd|^{-1}}^1 \frac 1 {t^3} dt \leq C |\log \Gd|^2,
\end{align*}
and by Lemma \ref{elem-propty} (iv)
\begin{align*}
|I_3(\Bx)|&\leq \left|\int_{p} ^{p_{N_1}}  \nabla \left(\frac  1  {|\Bx - (t,0,0)|} - \frac  1  {|\Bx + (t,0,0)|}  \right) ~\frac 1 {\sqrt {t^2- p^2}}~dt \right| \\
& \leq C \int_{p} ^{p_{N_1}} \frac 1 {\Gd^{5/4}} \frac 1 {\sqrt {t-p}}dt \leq C  \frac 1 {\Gd^{5/4} e^{ 1/(2\sqrt \Gd |\log \Gd|)}}\leq C.
\end{align*}
So far, we showed that
\beq\label{S1S3I1I3}
|S_1| + |S_3| + |I_1| + |I_3| \le C |\log \Gd|^2.
\eeq

We set
$$
\tilde{S}_2(\Bx)=\sum_{n=N_0}^{N_1-1}  \p_x\left(  \frac 1 {|\Bx - \Bp_n|}  - \frac 1 {|\Bx + \Bp_n|} \right)\frac{p_n-p_{n+1}}{\sqrt{p_n^2-p^2}},
$$
and shall prove
\beq\label{pro:result1}
\left|\tilde{S}_2(\Bx)- I_2(\Bx)\right| \leq  C ~\frac 1 {\sqrt \Gd+\Gr}.
\eeq

Let us first show that Lemma \ref{prop2} follows from \eqnref{S1S3I1I3} and \eqnref{pro:result1}. We observe from Lemma \ref{key_lemma} (i) that
$$
S_2(\Bx)=\tilde{S}_2(\Bx)\left(1+O(|\log \Gd|^{-1})\right).
$$
So we have
\begin{align*}
\p_x v & = S_1 + S_2 + S_3 \\
&= \tilde{S}_2 \left(1+O(|\log \Gd|^{-1})\right) + S_1 + S_3 \\
&= I_2 \left(1+O(|\log \Gd|^{-1})\right) + (\tilde{S}_2-I_2)\left(1+O(|\log \Gd|^{-1})\right) + S_1 + S_3 \\
&= \p_x v_0 \left(1+O(|\log \Gd|^{-1})\right) + R,
\end{align*}
where
$$
R = -(I_1+I_3) \left(1+O(|\log \Gd|^{-1})\right) + (\tilde{S}_2-I_2)\left(1+O(|\log \Gd|^{-1})\right) + S_1 + S_3.
$$
Since $\Gr \le |\log\Gd|^{-2}$, one can see from \eqnref{S1S3I1I3} and \eqnref{pro:result1} that
$$
|R| \le C \left( |\log\Gd|^2 + \frac 1 {\sqrt \Gd+\Gr} \right).
$$
Thanks to \eqnref{prop4_formul}, we now have
$$
\p_x v = \p_x v_0 \left(1+O(|\log \Gd|^{-1})\right) + R = \p_x v_0 \left(1+O(|\log \Gd|^{-1})\right),
$$
which we aim to prove.

The rest of this subsection is devoted to the proof of \eqnref{pro:result1}. For $N_0\leq n \leq N_1$, let
\begin{align*}
\gamma_n(\Bx) :&=\left|\p_x \left(\frac  1  {|\Bx - (p_n,0,0)|} - \frac  1  {|\Bx + (p_n,0,0)|}  \right)\frac{p_n-p_{n+1}}{\sqrt{p_n^2-p^2}}\right.\\
&\quad\left.- \int_{p_{n+1}}^{p_{n}} \p_x\left(\frac  1  {|\Bx - (t,0,0)|} - \frac  1  {|\Bx + (t,0,0)|}  \right) ~\frac 1 {\sqrt {t^2- p^2}}~dt\right|.
\end{align*}
Let
$$
f(t):= \p_x\left(\frac  1  {|\Bx - (t,0,0)|} - \frac  1  {|\Bx + (t,0,0)|}  \right) ~\frac 1 {\sqrt {t^2- p^2}}.
$$
By the mean-value property there is $t_n \in [p_{n+1},p_n]$ such that
$$
f(p_{n}) (p_n-p_{n+1})-\int_{p_{n+1}}^{p_{n}} f(t)dt =\frac{f'(t_n)}{2} (p_n-p_{n+1})^2.
$$
So we have
\begin{align*}
\gamma_n(\Bx) ~\leq &~   \frac{1}{2} \left| \p_t \p_{x} \Bigr(\frac  1  {|\Bx - (t,0,0)|} - \frac  1  {|\Bx + (t,0,0)|}   \Bigr) \Big|_{t=t_n} \right| \frac 1 {\sqrt {t_n^2 - p^2}} (p_n - p_{n+1})^2\\&
~+ \frac{1}{2} \left|\p_x \Bigr(\frac  1  {|\Bx - (t_n,0,0)|} - \frac  1  {|\Bx + (t_n,0,0)|}  \Bigr)\right| \frac {t_n } {( {t_n^2 - p^2})^{3/2}} (p_n - p_{n+1})^2\\
 := &~ \frac{1}{2} (\gamma_{n1}(\Bx) + \gamma_{n2}(\Bx)).
\end{align*}
Using \eqnref{xinRep}, one can show that
\beq
|\Bx\pm (t_n,0,0)|^2 \geq C(\Gr^2+ t_n^2), \quad\Bx\in R_\Gd
\eeq
for some constant independent of $n$. So we have
\beq\label{gamman1}
\Gg_{n1} \le \frac{C}{\Gr^3 + t_n^3} \frac{1}{\sqrt {t_n^2 - p^2}} (p_n - p_{n+1})^2
\eeq
and
\beq\label{gamman2}
\Gg_{n2} \le \frac{C}{\Gr^2 + t_n^2} \frac {t_n } {( {t_n^2 - p^2})^{3/2}} (p_n - p_{n+1})^2.
\eeq

If $n \leq N= [ \frac 1 {\sqrt {\Gd}}]$, then we have $t_n \approx 1/n$ and $|p_n - p_{n+1}| < C /n^2$ by Lemma \ref{lem:pq1overn}, and $p_n - p \ge C/n$ by Lemma \ref{elem-propty} (iii). So, we have
\begin{align*}
\sum_{n=N_0}^{N} \gamma_{n1} \leq  C \sum_{n=1} ^{N}   \frac 1 {(1/n^3 + \Gr^3)\sqrt {1/n^2}} ~ \frac 1 {n^4} \leq  C ~ \sum_{n=1}^{N} \frac 1  {1 + \Gr^3 n^3} .
\end{align*}
Note that if $\Gr \le \sqrt{\Gd}$, then
$$
\sum_{n=1}^{N} \frac 1  {1 + \Gr^3 n^3} \le N \le \frac{1}{\sqrt{\Gd}},
$$
while if $\Gr > \sqrt{\Gd}$, then
$$
\sum_{n=1}^{N} \frac 1  {1 + \Gr^3 n^3} \le \frac{1}{\Gr^3} \sum_{n=1}^{N} \frac{1}{n^3} \le \frac{C}{\Gr}.
$$
So, we have
$$
\sum_{n=N_0}^{N} \gamma_{n1} \leq \frac {C}{ \sqrt {\Gd} + \Gr}.
$$

If $N\leq n \leq N_1$, we have from \eqnref{pnp} that
$$
0 \le p_n - p_{n+1} = 2 p \frac {A^{-n-1} (1-A^{-1})} {(1-A^{-n-1}) (1-A^{-n-2})} \leq C_2  \Gd  A^{-n},
$$
and by Lemma \ref{elem-propty} (ii), $ p_n -p \geq C  \sqrt {\Gd} A^{-n}$.
Since $p_n >p$ for all $n$, we have
\begin{align*}
\sum_{n=N+1} ^{N_1-1} \gamma_{n1} &\leq C\sum_{n=N+1} ^{N_1-1}   \frac 1 {(p^{3} + \Gr^3) \sqrt {p_{n+1}^2 - p^2}} (p_n - p_{n+1})^2\notag\\
& \leq C\sum_{n=N+1}^{\infty}\frac 1 {(\Gd^{3/2} + \Gr^3) \Gd^{1/2} A^{-n/2}   } \Gd^2 A^{-2n}\notag\\
& \leq C\sum_{n=N+1} ^{\infty}\frac 1 {\Gd^{1/2} (\Gd^{3/2} + \Gr^3)   } \Gd^2 A^{-(3/2) n}\notag\\
& \leq C \frac {\Gd^{3/2}}{ \Gd^{3/2} + \Gr^3} \frac 1 {\sqrt{\Gd}} \leq \frac {C}{ \sqrt {\Gd} + \Gr}.
\end{align*}
Thus, we have
$$
\sum_{n=N_0} ^{N_1} \gamma_{n1} \leq  \frac C {\Gr + \sqrt \Gd}.
$$

Similarly one can show that
$$
\sum_{n=N_0} ^{N_1} \gamma_{n2} \leq  \frac C {\Gr + \sqrt \Gd}.
$$
This completes the proof of \eqnref{pro:result1}. \qed

\section{Proofs of Theorem \ref{main} and Theorem \ref{main2}}\label{sec:proof}

{\sl Proof of Theorem \ref{main}}. It is helpful to recall that $\Ge=2\Gd$. We have from \eqnref{old_result} and \eqnref{hconst} that
\begin{align}
u &= \frac {u |_{\p D_1}- u |_{\p D_2} }{h |_{\p D_1}- h |_{\p D_2}}  h +  g \nonumber \\
& = \frac{1}{2} ( u |_{\p D_2}- u |_{\p D_1} ) (4\pi\sum_{n=0}^{\infty} q_n) h + g \nonumber \\
&= \frac{1}{2} ( u |_{\p D_2}- u |_{\p D_1} ) v + g. \label{511}
\end{align}
We emphasize that $|\nabla g|$ is bounded on any bounded subset of $\Rbb^3 \setminus (D_1 \cup D_2)$ regardless of $\Gd$ as explained in Introduction. Since $h$ is constant on $\p D_1$ and $\p D_2$, one can see from \eqnref{h:basic} and \eqnref{h:sameradii} that
\begin{align*}
u|_{\p D_2}- u|_{\p D_1} & = -\int_{\p (D_1 \cup D_2)}  H  \frac{\p h}{\p \nu} d\Gs \nonumber \\
&= -\int_{\p (D_1 \cup D_2)}  H  \frac{\p h}{\p \nu} d\Gs + \int_{\p (D_1 \cup D_2)}  \frac{\p H}{\p \nu} h d\Gs \nonumber \\
&= -\frac{1}{\sum_{n=0}^{\infty} q_n} \sum_{n=0}^{\infty} q_n \int_{\p (D_1 \cup D_2)} \Big[ H(\Bx) \frac{\p}{\p \nu}\left(\GG(\Bx - \Bp_n) - \GG(\Bx + \Bp_n) \right)  \nonumber \\
& \qquad\qquad  - \frac{\p H}{\p \nu} \left(\GG(\Bx - \Bp_n) - \GG(\Bx + \Bp_n) \right) \Big] d\Gs \nonumber \\
& =\frac{C_H^\Ge}{2\sum_{n=0}^{\infty}  q_n}.
\end{align*}
It then follows from \eqnref{q:sum} that
\beq\label{udiff:Cep}
u |_{\p D_2}- u |_{\p D_1}=\frac {C_H^\Ge} { |\log \Gd|}\left(1+ O(|\log \Gd|^{-1})\right),
\eeq
where $O(|\log \Gd|^{-1})$ is independent of $H$. So we obtain from \eqnref{511}
$$
\nabla u =  \frac {C_H^\Ge} {2|\log \Gd|} \left(1+ O(|\log \Gd|^{-1})\right) \nabla v +\nabla g ,
$$
and from \eqnref{asymp_sing}
\beq\label{inter}
\nabla u(\Bx) =  \frac{C_H^\Ge} { |\log \Gd| (2\Gd + \Gr(\Bx)^2)} \left((1,0,0) + O(|\log \Gd|^{-1}) \right) +\nabla g (\Bx),
\eeq
and hence \eqnref{asymp1} is proved. This completes the proof. \qed

\medskip
\noindent{\sl Proof of Theorem \ref{main2}}. By \eqnref{pq1overn}, we have for $n \le N=N(\Gd)$
\begin{align*}
& \left| q_n (H(\Bp_n)-H(-\Bp_n)) - \frac{1}{n+1} \left( H\left(\frac{1}{n+1},0,0\right) - H\left(-\frac{1}{n+1},0,0\right) \right) \right| \\
& \le \left| q_n - \frac{1}{n+1} \right| \left| H(\Bp_n)-H(-\Bp_n) \right| \\
& \quad + \frac{1}{n+1} \left( \left| H(\Bp_n) - H\left(\frac{1}{n+1},0,0\right) \right| +
\left| H(-\Bp_n) - H\left(-\frac{1}{n+1},0,0\right) \right| \right) \\
& \le C\sqrt{\Gd} \left( \sqrt{\Gd} + \frac{1}{n+1} \right) + C \sqrt{\Gd} \frac{1}{n+1}.
\end{align*}
So we have
$$
\left| \sum_{n=0}^{N-1} q_n (H(\Bp_n)-H(-\Bp_n)) - \sum_{n=1}^{N} \frac{1}{n} \left( H\left(\frac{1}{n},0,0\right) - H\left(-\frac{1}{n},0,0\right) \right) \right| \le C\sqrt{\Gd} |\log\Gd|.
$$

On the other hand, since $p_n$ is decreasing, it follows from \eqnref{q:sum} that
$$
\left| \sum_{n=N}^{\infty} q_n (H(\Bp_n)-H(-\Bp_n)) \right| \le C p_N \sum_{n=N}^{\infty} q_n \le C\sqrt{\Gd} |\log\Gd|.
$$
We also have
$$
\left| \sum_{n=N+1}^{\infty} \frac{1}{n} \left( H\left(\frac{1}{n},0,0\right) - H\left(-\frac{1}{n},0,0\right) \right) \right| \le C \sum_{n=N+1}^{\infty} \frac{1}{n^2} \le C \sqrt{\Gd}.
$$
Combining above estimates, we obtain \eqnref{0.3}.

The formula \eqnref{intensity} is an immediate consequence of \eqnref{asymp1} and \eqnref{0.3}. This completes the proof. \qed


\end{document}